\documentclass[a4paper,12pt]{amsart}
\usepackage{amssymb, amsmath, graphicx}
\usepackage[curve]{xypic}

\providecommand{\Ker}{\textnormal{Ker}}
\providecommand{\IIm}{\textnormal{Im}}

\providecommand{\Cyl}{\textnormal{Cyl}}

\providecommand{\Ker}{\textnormal{Ker}}

\providecommand{\ch}{\textnormal{ch}}

\providecommand{\fl}{\textnormal{fl}}

\begin{document}

\title{Relative differential cohomology}
\author{Fabio Ferrari Ruffino}
\address{Departamento de Matem\'atica - Universidade Federal de S\~ao Carlos - Rod.\ Washington Lu\'is, Km 235 - C.P.\ 676 - 13565-905 S\~ao Carlos, SP, Brasil}
\email{ferrariruffino@gmail.com}

\begin{abstract}
Let $h^{\bullet}$ be a cohomology theory and $\hat{h}^{\bullet}$ the natural differential refinement, as defined by Hopkins and Singer. We consider the possible definitions of the relative differential cohomology groups, generalizing the analogous picture for the Deligne cohomology, and we show the corresponding long exact sequence in each case. \\ $ $ \\
\textsc{MSC2010: 55N20, 53C08}
\end{abstract}

\maketitle

\newtheorem{Theorem}{Theorem}[section]
\newtheorem{Lemma}[Theorem]{Lemma}
\newtheorem{Corollary}[Theorem]{Corollary}
\newtheorem{Rmk}[Theorem]{Remark}
\newtheorem{Def}{Definition}[section]
\newtheorem{ThmDef}[Theorem]{Theorem - Definition}

\section{Introduction}

Let $h^{\bullet}$ be a cohomology theory and $\hat{h}^{\bullet}$ the natural differential refinement, as defined in \cite{HS, Upmeier, UpmeierThesis}. In \cite{Upmeier, UpmeierThesis} a relative version of $\hat{h}^{\bullet}$ has been defined, considering classes on a pair of manifolds $(X, A)$ whose restriction to $A$ is trivial: nevertheless, when considering the differential refinement of ordinary cohomology, this version is not the only interesting one. It is more natural to consider the case in which the restriction to $A$ is only topologically trivial, so that the curvature is exact but not necessarily vanishing, as in the first definition of relative Cheeger-Simons character given in \cite{BT} (see also \cite{BB, Becker}). For example, the $B$-field in string theory can be described as a Deligne cohomology class of degree $2$ in the space-time, which is topologically trivial on a spin$^{c}$ D-brane world-volume, because of the Freed-Witten anomaly \cite{FW, BFS, FR4}. In this case, imposing that the $B$-field is vanishing on the world-volume is an unnecessary restriction. Moreover, in \cite{BT} the authors showed that it is possible to define the relative Cheeger-Simons characters in such a way that they fit into a long exact sequence completely made by differential cohomology groups, and it is natural to inquire if such a definition can be generalized.

We have shown in \cite{FR3} that there are actually four possible inequivalent ways to define the relative Deligne cohomology groups, two of which being the most meaningful and corresponding to the two kinds of differential characters defined in \cite{BT, Becker}. We generalize this picture to any cohomology theory and we show in each case the corresponding long exact sequence. Such a generalization can be applied for example to describe the B-field when it is thought of as an element of a generalized cohomology theory different from the ordinary one \cite{DFM}, or any background field that must be topologically trivial on a subspace in order to delete an anomaly.

The paper is organized as follows. In section \ref{RelDiffCoh} we state the possible definitions of relative differential cohomology for a pair of manifolds and we show the corresponding long exact sequences. In section \ref{ExSequences} we construct the Bockstein map of each exact sequence and we prove the exactness. In section \ref{GenericMap} we generalize the previous constructions to any map of manifolds, not necessarily a closed embedding.

\section{Relative differential cohomology}\label{RelDiffCoh}

We fix a cohomology theory $h^{\bullet}$, represented by an $\Omega$-spectrum $(E_{n}, e_{n}, \varepsilon_{n})$, where $e_{n}$ is the marked point of $E_{n}$ and $\varepsilon_{n}: (\Sigma E_{n}, \Sigma e_{n}) \rightarrow (E_{n+1}, e_{n+1})$ is the structure map, whose adjoint $\tilde{\varepsilon}_{n}: E_{n} \rightarrow \Omega_{e_{n+1}} E_{n+1}$ is a homeomorphism. We also fix real singular cocycles $\iota_{n} \in C^{n}(E_{n}, e_{n}, \mathfrak{h}^{\bullet}_{\mathbb{R}})$ representing the Chern character of $h^{\bullet}$, such that $\iota_{n-1} = \int_{S^{1}} \varepsilon_{n}^{*}\iota_{n}$ \cite{Upmeier}. Moreover, for $\{*\}$ a space with one point, we call $\mathfrak{h}^{\bullet}_{\mathbb{R}} := h^{\bullet}(\{*\}) \otimes_{\mathbb{Z}} \mathbb{R}$. In this section we state the possible definitions of the relative differential cohomology groups and we show the corresponding long exact sequences, postponing the construction of the Bockstein maps and the proofs to the next section.

\subsection{Relative differential cohomology of type I} This definition is the one considered in \cite{Upmeier}, and generalizes to any cohomology theory the first one of \cite{FR3}. In the following, we call \emph{pair of smooth manifolds} a pair $(X, A)$ such that $X$ is a smooth manifold and $A \subset X$ is a \emph{closed} embedded submanifold; $X$ is allowed to have a boundary, in which case $A$ is a neat submanifold.
\begin{Def}\label{DiffFunct} If $(X, A)$ is a pair of smooth manifolds, $(Y, y_{0})$ a topological space with marked point, $V^{\bullet}$ a graded real vector space and $\kappa_{n} \in C^{n}(Y, y_{0}, V^{\bullet})$ a real singular cocycle, a \emph{relative differential function of type I} from $(X, A)$ to $(Y, y_{0}, \kappa_{n})$ is a triple $(f, h, \omega)$ such that:
\begin{itemize}
	\item $f: (X, A) \rightarrow (Y, y_{0})$ is a continuous function;
	\item $h \in C^{n-1}_{sm}(X, A; V^{\bullet})$ (`sm' means smooth);
	\item $\omega \in \Omega^{n}_{cl}(X, A; V^{\bullet})$ (i.e., $\omega\vert_{A} = 0$)
\end{itemize}
satisfying, for $\chi: \Omega^{\bullet}(X, A; V^{\bullet}) \rightarrow C^{\bullet}(X, A; V^{\bullet})$ the natural homomorphism:
\begin{equation}\label{DiffFunctEq}
	\delta^{n-1}h = \chi^{n}(\omega) - f^{*}\kappa_{n}.
\end{equation}
Moreover, a \emph{homotopy between two relative differential functions} $(f_{0}, h_{0}, \omega)$ and $(f_{1}, h_{1}, \omega)$ is a relative differential function $(F, H, \pi^{*}\omega): (X \times I, A \times I) \rightarrow (Y, y_{0}, \kappa_{n})$, such that $F$ is a homotopy between $f_{0}$ and $f_{1}$, $H\vert_{(X \times \{i\}, A \times \{i\})} = h_{i}$ for $i = 0,1$, and $\pi: X \times I \rightarrow X$ is the natural projection.
\end{Def}
A cochain $h \in C^{n-1}(X, A)$ is uniquely represented by a cochain $h \in C^{n-1}(X)$ such that $h\vert_{A} = 0$. Hence, in the previous definition, we could equivalently require that $(f, h, \omega)$ is a differential function from $X$ to $(Y, \kappa_{n})$ such that $(f, h, \omega)\vert_{A} = (c_{y_{0}}, 0, 0)$, being $c_{y_{0}}$ the constant map with value $y_{0}$. It follows from the definition that a homotopy is constant on $A \times I$, i.e.\ $(F, H, \pi^{*}\omega)\vert_{A \times \{t\}} = (c_{e_{n}}, 0, 0)$ for any $t \in I$. We define the group $\hat{h}_{I}^{n}(X, A)$ in the following way: as a set, $\hat{h}_{I}^{n}(X, A)$ contains the homotopy classes of type I relative differential functions $(f, h, \omega): (X, A) \rightarrow (E_{n}, e_{n}, \iota_{n})$; the abelian group structure is defined similarly to the absolute case \cite{UpmeierThesis, Upmeier}. The corresponding long exact sequence is the following:
\begin{equation}\label{LongExact1}
\begin{split}
	\cdots & \longrightarrow \hat{h}^{p-1}_{\fl}(X, A) \longrightarrow \hat{h}^{p-1}_{\fl}(X) \longrightarrow \hat{h}^{p-1}_{\fl}(A) \\
	& \longrightarrow \hat{h}^{p}_{I}(X, A) \longrightarrow \hat{h}^{p}(X) \longrightarrow \hat{h}^{p}(A) \\
	& \longrightarrow h^{p+2}(X, A) \longrightarrow h^{p+2}(X) \longrightarrow h^{p+2}(A) \longrightarrow \cdots.
\end{split}
\end{equation}

\subsection{Relative differential cohomology of type II} This definition generalizes to any cohomology theory the second one of \cite{FR3}. It is the most interesting for applications. In the groups of type I we require that the differential function is vanishing on $A$. For Deligne cohomology, in the groups of type II the differential class on $A$ must be topologically trivial, and a specific trivialization is chosen as a part of the data, therefore the class on $A$ is represented by a Deligne cocycle of the form $(1, 0, \ldots, 0, \rho)$. Here we generalize this requirement, replacing the representatives of the form $(1, 0, \ldots, 0, \rho)$ with the ones of the form $(c_{e_{n}}, \chi(\rho), d\rho)$. Moreover, as for type I, a homotopy must be constant on $A$: this means that it must be of the form $(c_{e_{n}}, \pi^{*}\chi(\rho), \pi^{*}d\rho)$ on $A \times I$, where $\pi: A \times I \rightarrow A$ is the projection. Hence, with the same data of definition \ref{DiffFunct}, we define a \emph{relative differential function of type II} from $(X, A)$ to $(Y, y_{0}, \kappa_{n})$ as a quadruple $(f, h, \omega, \rho)$ such that:
\begin{itemize}
	\item $(f, h, \omega)$ is a differential function from $X$ to $(Y, \kappa_{n})$;
	\item $\rho \in \Omega^{n-1}(A; V^{\bullet})$;
	\item $(f, h, \omega)\vert_{A} = (c_{y_{0}}, \chi(\rho), d\rho)$.
\end{itemize}
A homotopy between $(f_{0}, h_{0}, \omega, \rho)$ and $(f_{1}, h_{1}, \omega, \rho)$ is a differential function $(F, H, \pi^{*}\omega, \pi^{*}\rho): (X \times I, A \times I) \rightarrow (Y, y_{0}, \kappa_{n})$, such that $F$ is a homotopy between $f_{0}$ and $f_{1}$, $H\vert_{(X \times \{i\}, A \times \{i\})} = h_{i}$ for $i = 0,1$, and $\pi: X \times I \rightarrow X$ is the natural projection. Considering homotopy classes of type II differential functions, instead of those of type I, we define the group $\hat{h}_{II}^{n}(X, A)$. The corresponding long exact sequence is the following:
\begin{equation}\label{LongExact2}
\begin{split}
	\cdots & \longrightarrow \hat{h}^{p-1}_{\fl}(X, A) \longrightarrow \hat{h}^{p-1}_{\fl}(X) \longrightarrow \hat{h}^{p-1}(A) \\
	& \longrightarrow \hat{h}^{p}_{II}(X, A) \longrightarrow \hat{h}^{p}(X) \longrightarrow h^{p+1}(A) \\
	& \longrightarrow h^{p+2}(X, A) \longrightarrow h^{p+2}(X) \longrightarrow h^{p+2}(A) \longrightarrow \cdots.
\end{split}
\end{equation}
There is a natural embedding
\begin{equation}\label{MapPP1h}
	\varphi_{I, II}: \hat{h}^{p}_{I}(X, A) \hookrightarrow \hat{h}^{p}_{II}(X, A),
\end{equation}
which actually extends to a morphism of exact sequences from \eqref{LongExact1} to \eqref{LongExact2}.

\subsection{Relative differential cohomology of type III} This definition is less interesting than the previous ones, but we put it anyway for completeness. We must require that a relative class is topologically trivial on $A$, with no reference to any specific trivialization of the differential class. Hence, with the same data of definition \ref{DiffFunct}, we define a \emph{relative differential function of type III} from $(X, A)$ to $(Y, y_{0}, \kappa_{n})$ as a differential function from $X$ to $(Y, \kappa_{n})$ such that $f\vert_{A} = c_{y_{0}}$. Homotopies are defined as in the previous cases. The corresponding long exact sequence is the following:
\begin{equation}\label{LongExact3}
\begin{split}
	\cdots & \longrightarrow \hat{h}^{p-1}_{\fl}(X, A) \longrightarrow \hat{h}^{p-1}_{\fl}(X) \longrightarrow h^{p}(A) \\
	& \longrightarrow \hat{h}^{p}_{III}(X, A) \longrightarrow \hat{h}^{p}(X) \longrightarrow h^{p+1}(A) \\
	& \longrightarrow h^{p+2}(X, A) \longrightarrow h^{p+2}(X) \longrightarrow h^{p+2}(A) \longrightarrow \cdots.
\end{split}
\end{equation}
There is a natural surjective map
\begin{equation}\label{MapP10h}
	\varphi_{II, III}: \hat{h}^{p}_{II}(X, A) \rightarrow \hat{h}^{p}_{III}(X, A),
\end{equation}
which actually extends to a morphism of exact sequences from \eqref{LongExact2} to \eqref{LongExact3}.

\subsection{Relative differential cohomology of type IV} In the group $C^{\bullet}(X, \mathbb{R})$ we consider the subgroup $Z^{\bullet}_{\ch}(X, \mathbb{R})$, containing the cocycles that represent a class belonging to the image of the Chern character of the theory $h^{\bullet}$ (for ordinary cohomology such a subgroup is generated by $Z^{\bullet}(X, \mathbb{Z})$ and the real coboundaries). We can consider a definition of differential function analogous to the usual one, with the only difference that, in the triple $(f, h, \omega)$, we think of $h$ as an element of $C^{n-1}(X, \mathbb{R})/Z^{n-1}_{\ch}(X, \mathbb{R})$ instead of $C^{n-1}(X, \mathbb{R})$. This has no effects on differential cohomology, since a class $[(c_{e_{n}}, h, 0)]$ is vanishing if and only if $h \in Z^{\bullet}_{\ch}(X, \mathbb{R})$. Nevertheless, this quotient has some effects on the relative groups. In particular, we define a \emph{relative differential function of type IV} replacing $h\vert_{A} = 0$ in the type I definition (v.\ comments after definition \ref{DiffFunct}) with $[h]\vert_{A} = 0$ as an element of $C^{n-1}(X, \mathbb{R})/Z^{n-1}_{\ch}(X, \mathbb{R})$. With this definition we get a long exact sequence all made by differential cohomology groups, while the ones we obtained up to now contain some topological groups and some differential extensions. The sequence is the following:
\begin{equation}\label{LongExact4}
\begin{split}
	\cdots & \longrightarrow \hat{h}_{VI}^{p-1}(X, A) \longrightarrow \hat{h}^{p-1}(X) \longrightarrow \hat{h}^{p-1}(A) \\
	& \longrightarrow \hat{h}^{p}_{IV}(X, A) \longrightarrow \hat{h}^{p}(X) \longrightarrow \hat{h}^{p}(A) \\
	& \longrightarrow \hat{h}_{IV}^{p+1}(X, A) \longrightarrow \hat{h}^{p+1}(X) \longrightarrow \hat{h}^{p+1}(A) \longrightarrow \cdots.
\end{split}
\end{equation}

\section{The exact sequences}\label{ExSequences}

We now construct the Bockstein maps of the long exact sequences stated in the previous section, and we prove the exactness in each case.

\subsection{Homotopy extension property}\label{HomExtProp} In the following we will need the homotopy extension property for differential functions. We call $\Cyl \, A := A \times I$ and, for a pair $(X,A)$, we call $\Cyl(X,A)$ the union $X \cup \Cyl \, A$ identifying $A \subset X$ with $A \times \{0\} \subset \Cyl \, A$. In general $\Cyl(X,A)$ is not a manifold, nevertheless we will deal with differential functions $(f, h, \omega): \Cyl(X,A) \rightarrow (E_{n}, \iota_{n})$, defined in the following way:
\begin{itemize}
	\item $f: \Cyl(X,A) \rightarrow E_{n}$ is a continuous function.
	\item $\omega \in \Omega^{n}_{cl}(X; \mathfrak{h}^{\bullet}_{\mathbb{R}})$, and it defines a smooth cocycle $\chi^{n}(\omega)$ on $\Cyl(X, A)$ as follows. $\Cyl(X, A)$ is a subspace of $X \times I$. Let us consider the pull-back $\pi_{X}^{*}\omega$ on $X \times I$ and the embedding $\iota: \Cyl(X, A) \hookrightarrow X \times I$. A simplex $\sigma: \Delta^{n} \rightarrow \Cyl(X, A)$ is defined to be smooth if and only if the composition $\iota \circ \sigma: \Delta^{n} \rightarrow X \times I$ is. The smooth cochain $\chi^{n}(\omega)$ on $\Cyl(X, A)$ is defined by $\chi^{n}(\omega)(\sigma) := \chi^{n}(\pi_{X}^{*}\omega)(\iota \circ \sigma)$ (i.e., up to the embedding $\iota$, $\chi^{n}(\omega)$ is the restriction of $\chi^{n}(\pi_{X}^{*}\omega)$ to the smooth chains whose image is contained in $\Cyl(X, A)$).
	\item $h \in C^{n-1}_{sm}(\Cyl(X,A); \mathfrak{h}^{\bullet}_{\mathbb{R}})$ and it satisfies $\delta^{n-1}h = \chi^{n}(\omega) - f^{*}\iota_{n}$.
\end{itemize}
The same considerations apply if we iterate the cylinder, for example considering the space $\Cyl(\Cyl(X,A),\Cyl\,A)$, and so on. We call ``manifold with cylinders'' a space obtained in this way. Let $f: X \rightarrow Y$ be a function between manifolds with cylinders. By definition $X$ and $Y$ are contained in manifolds of the form $X' \times I^{k}$and $Y' \times I^{h}$. We call $f$ smooth if and only if the composition $\iota \circ f: X \rightarrow Y' \times I^{h}$ is the restriction of a smooth function defined on a neighborhood of $X$ in $X' \times I^{k}$.
\begin{Lemma}\label{HEP} Given a pair of smooth manifolds with cylinders $(X, A)$, the homotopy extension property holds for differential functions, i.e., given $(f, h, \omega): \Cyl(X,A) \rightarrow (Y, \kappa_{n})$, there exists a homotopy $(F, H, \pi_{X}^{*}\omega): X \times I \rightarrow (Y, \kappa_{n})$ extending $(f, h, \omega)$.
\end{Lemma}
\paragraph{Proof:} Since a pair of smooth manifolds with cylinders is a CW-pair, the topological homotopy extension property holds, therefore we can extend $f$ to $F: X \times I \rightarrow Y$. Since $F$ is homotopic to $f \circ \pi_{X}$, we have $[\chi^{n}(\pi_{X}^{*}\omega) - F^{*}\iota_{n}] = \pi_{X}^{*}[\chi^{n}(\omega) - f^{*}\iota_{n}] = [\pi_{X}^{*}\delta (h\vert_{X})] = 0$, hence there exists $H' \in C^{n-1}_{sm}(X \times I; \mathfrak{h}^{\bullet}_{\mathbb{R}})$ such that $\delta H' = \chi^{n}(\pi^{*}\omega) - F^{*}\iota_{n}$. The homotopy extension property is equivalent to $\Cyl(X,A)$ being a retract of $X \times I$. Moreover, since $\Cyl(X,A)$ is a closed subspace, a continuous retraction $r: X \times I \rightarrow \Cyl(X,A)$ can be turned into a smooth one, so that the pull-back of smooth cochains is well-defined. Let $K := r^{*}(h - H'\vert_{X \cup \Cyl \, A})$: one has $\delta K = 0$ by construction. We define $H := H' + K$. Then $\delta H = \chi^{n}(\pi^{*}\omega) - F^{*}\iota_{n}$ and $H\vert_{\Cyl(X,A)} = h$. $\square$ \\

\subsection{Type I}\label{TypeISubS} The exactness of \eqref{LongExact1} has been already proven in \cite{Upmeier}. In particular, the relative flat theory $\hat{h}^{p}_{\fl}(X, A)$ is by definition the subgroup of $\hat{h}^{p}_{I}(X, A)$ made by classes with vanishing curvature, and, being a cohomology theory, it defines a long exact sequence. Hence the Bockstein map of \eqref{LongExact1} is defined as the composition of the Bockstein map of the flat theory $\hat{h}^{p-1}_{\fl}(A) \rightarrow \hat{h}^{p}_{\fl}(X, A)$ with the immersion $\hat{h}^{p}_{\fl}(X, A) \hookrightarrow \hat{h}^{p}_{I}(X, A)$.

We define the relative groups of type I in an alternative way, which will be suitable to be generalized to type II in order to construct the Bockstein map of \eqref{LongExact2}. When we consider the relative Deligne cohomology, a relative class of this type is represented by a cocycle on $X$ together with a geometrical trivialization of its restriction to $A$ \cite{FR3}. Here we can repeat an analogous construction, considering the trivialization as a homotopy of differential functions.
\begin{Def} Given a differential function $(f, h, \omega): X \rightarrow (E_{n}, \iota_{n})$, a \emph{geometric trivialization} of $(f, h, \omega)$ is a homotopy $(F, H, \pi^{*}\omega): X \times I \rightarrow (E_{n}, \iota_{n})$ between $(f, h, \omega)$ and the trivial function $(c_{e_{n}}, 0, 0)$.
\end{Def}
It follows from the definition that $\omega = 0$.
\begin{Def}\label{GroupILinha} The group $\hat{h}^{n}_{I'}(X, A)$ contains the homotopy classes of differential functions $(f, h, \omega): \Cyl(X,A) \rightarrow (E_{n}, \iota_{n})$ such that $(f, h, \omega)\vert_{\Cyl A}$ is a geometric trivialization of $(f, h, \omega)\vert_{A}$. A homotopy $(F, H, \pi^{*}\omega): \Cyl(X, A) \times I \rightarrow (E_{n}, \iota_{n})$ between two such functions is required to satisfy $(F, H, \pi^{*}\omega)\vert_{A \times \{1\} \times I} = (c_{e_{n}}, 0, 0)$.
\end{Def}
It follows that $\omega\vert_{A} = 0$. We can show that $\hat{h}^{n}_{I'}(X, A) \simeq \hat{h}_{I}^{n}(X, A)$ canonically. There is a natural morphism $\varphi: \hat{h}_{I}^{n}(X, A) \rightarrow \hat{h}^{n}_{I'}(X, A)$: given $[(f, h, \omega)] \in \hat{h}_{I}^{n}(X, A)$, since by definition $(f, h, \omega)\vert_{A} = (c_{e_{n}}, 0, 0)$, we extend $(f, h, \omega)$ to $\Cyl(X,A)$ in such a way that $(f, h, \omega)\vert_{\Cyl\,A} = (c_{e_{n}}, 0, 0)$ (v.\ \cite{Upmeier} lemma 2.2). It is easy to verify that $\varphi$ is well-defined up to homotopy and that respects the sum.
\begin{Theorem}\label{IsoI} The morphism $\varphi: \hat{h}_{I}^{n}(X, A) \rightarrow \hat{h}^{n}_{I'}(X, A)$ is an isomorphism.
\end{Theorem}
\paragraph{Proof:} \subparagraph{\emph{Injectivity}} Let $\varphi([(f, h, \omega)]) = 0$. Then there exists a homotopy $(F, H, \pi^{*}\omega)$ between $(f, h, \omega)$, extended to $\Cyl(X,A)$, and $(c_{e_{n}}, 0, 0)$, such that $(F, H, \pi^{*}\omega)\vert_{A \times \{1\} \times I} = (c_{e_{n}}, 0, 0)$. We now apply the homotopy extension property extending $(F, H, \pi^{*}\omega)$ to a function $(F', H', \pi^{*}\omega): X \times I \times I \rightarrow (E_{n}, \iota_{n})$ such that $(F', H', \pi^{*}\omega)\vert_{A \times I \times \{0\}}$, $(F', H', \pi^{*}\omega)\vert_{A \times \{1\} \times I}$ and $(F', H', \pi^{*}\omega)\vert_{A \times I \times \{1\}}$ are all equal to $(c_{e_{n}}, 0, 0)$. Therefore, composing the homotopies $(F', H', \pi^{*}\omega)\vert_{X \times I \times \{0\}}$, $(F', H', \pi^{*}\omega)\vert_{X \times \{1\} \times I}$ and $(F', H', \pi^{*}\omega)\vert_{X \times I \times \{1\}}$ we get a homotopy between $(f, h, \omega)$ and $(c_{e_{n}}, 0, 0)$ which is trivial on $A$, thus we get a homotopy of relative differential functions between $(f, h, \omega): (X, A) \rightarrow (E_{n}, e_{n}, \iota_{n})$ and $(c_{e_{n}}, 0, 0)$: this proves that $[(f, h, \omega)] = 0$ in $\hat{h}_{I}^{n}(X, A)$.
\subparagraph{\emph{Surjectivity}} Let us consider $[(f, h, \omega)] \in \hat{h}^{n}_{I'}(X, A)$, in particular $(f, h, \omega): \Cyl(X, A) \rightarrow (E_{n}, \iota_{n})$ and $(f, h, \omega)\vert_{A \times \{1\}} = (c_{e_{n}}, 0 ,0)$. We construct a homotopy $(F, H, 0): A \times I \times I \rightarrow (E_{n}, \iota_{n})$ between $(f, h, \omega)\vert_{A \times I}$ and $(c_{e_{n}}, 0, 0)$, such that $(F, H, 0)\vert_{A \times \{1\} \times I} = (c_{e_{n}}, 0, 0)$. We define:
	\[F(a, u, t) := \left\{ \begin{array}{lll} f(a, u+t) & & u + t \leq 1 \\ f(a, 1) & & u + t \geq 1. \end{array} \right.
\]
In this way $F(a, u, 0) = f(a, u)$ and $F(a, u, 1) = F(a, 1, t) = f(a, 1) = e_{n}$. We call $\pi_{A \times I \times \{0\}}: A \times I \times I \rightarrow A \times I \times \{0\}$ the projection, and similarly for the other cases. We have that $F$ is homotopic to $f \circ \pi_{A \times I \times \{0\}}$, therefore $[F^{*}\iota_{n}] = -\pi_{A \times I \times \{0\}}^{*}[\delta^{n-1}h] = 0$, hence there exists $H'' \in C^{n-1}_{sm}(A \times I \times I, \mathfrak{h}^{\bullet}_{\mathbb{R}})$ such that $\delta H'' = -F^{*}\iota_{n}$. We must now replace $H''$ by a cocycle $H$ such that $H\vert_{A \times I \times \{0\}} = h$ and $H\vert_{A \times I \times \{1\}} = H\vert_{A \times \{1\} \times I} = 0$. We do it in three steps:
\begin{itemize}
	\item We call $K'' := \pi_{A \times \{1\} \times I}^{*}(H''\vert_{A \times \{1\} \times I})$. Then $\delta K'' = 0$, hence we define $H' = H'' - K''$. In this way $\delta H' = -F^{*}\iota_{n}$ and $H'\vert_{A \times \{1\} \times I} = 0$.
	\item We call $K' := \pi_{A \times I \times \{0\}}^{*}(h - H'\vert_{A \times I \times \{0\}})$: then $\delta K' = 0$. For $H := H' + K'$ we get $\delta H = -F^{*}\iota_{n}$ and $H\vert_{A \times I \times \{0\}} = h$, keeping $H\vert_{A \times \{1\} \times I} = 0$.
	\item Up to now $H\vert_{A \times I \times \{1\}} = h'$ with $\delta h' = 0$ and $h'\vert_{A \times \{1\} \times \{1\}} = 0$. Since $A \times \{1\} \times \{1\}$ is a deformation retract of $A \times I \times \{1\}$, the fact that $h'\vert_{A \times \{1\} \times \{1\}} = 0$ implies that the cohomology class of $h'$ vanishes on the whole $A \times I \times \{1\}$, i.e.\ $h' = \delta k'$. By \cite[Lemma 2.3]{Upmeier}, we can choose a homotopy between $(c_{e_{n}}, \delta k', 0)$ and $(c_{e_{n}}, 0, 0)$ of the form $(c_{e_{n}}, K', 0)$, with $K'\vert_{A \times \{1\} \times I} = 0$. We compose $(F, H, 0)$ with this homotopy, and for simplicity we still call the result $(F, H, 0)$. Now $H\vert_{A \times I \times \{1\}} = 0$.
\end{itemize}
Thanks to the homotopy extension property of the pair $(\Cyl(X, A), \Cyl\,A)$, we extend $(F, H, 0)$ to $(F, H, \pi^{*}\omega): \Cyl(X,A) \times I \rightarrow (E_{n}, \iota_{n})$. By construction the function $(F, H, \pi^{*}\omega): \Cyl(X,A) \times \{1\} \rightarrow (E_{n}, \iota_{n})$ represents a class lying in the image of $\varphi$ that, considering definition \ref{GroupILinha}, is homotopic to $(f, h, \omega)$: it follows that the class $[(f, h, \omega)]$ we started from belongs to the image of $\varphi$. $\square$ \\

\paragraph{\textbf{Remark:}} Given a class $[(f, h, \omega)] \in \hat{h}^{n}_{I'}(X, A)$, by the homotopy extension property applied to $(X, A)$ we can extend $(f, h, \omega)$ to $(\Phi, \mathcal{H}, \pi_{X}^{*}\omega): X \times I \rightarrow (E_{n}, \iota_{n})$. The restriction to $X \times \{1\}$ is a class $[(g, h, \omega)] \in \hat{h}^{n}_{I}(X, A)$ by construction. We claim that $\varphi[(g, h, \omega)] = [(f, h, \omega)]$. In fact, in the proof of surjectivity, we can extend $(F, H, 0)$ to the whole $\Cyl\,X \times I$, applying the homotopy extension property to the pair $(\Cyl\,X, \Cyl\,A)$ instead of $(\Cyl(X, A), \Cyl\,A)$. It follows that the homotopies $(F, H, \pi^{*}\omega)\vert_{X \times \{1\} \times I}$ and $(F, H, \pi^{*}\omega)\vert_{X \times I \times \{1\}}$ compose to a homotopy (with respect to def.\ \ref{DiffFunct}) between $(g, h, \omega)$ and $(F, H, \pi^{*}\omega)\vert_{X \times \{0\} \times \{1\}}$, the latter being by construction a representative of $\varphi^{-1}[(f, h, \omega)]$. $\square$ \\

The Bockstein map of \eqref{LongExact1} can be visualized in this way: given $\alpha \in \hat{h}_{\fl}^{n-1}(A, x_{0})$, we consider the suspension isomorphism in the flat theory and we get $\tilde{\alpha} \in \hat{h}_{\fl}^{n}(SA, *)$. For $p: A \times I \rightarrow SA$ the natural projection, such that $p(A \times \{0, 1\}) = \{*\}$, we consider $p^{*}\tilde{\alpha}$ and we extend it trivially on $\Cyl(X, A)$. The class we get is the image of $\alpha$ in $\hat{h}^{n}_{I'}(X, A)$. Thanks to the previous theorem, we get a corresponding class in $\hat{h}^{n}_{I}(X, A)$.

\subsection{Type II} We now adapt the construction of the previous paragraph to the groups of type II, in order to define the Bockstein map of \eqref{LongExact2} and show the exactness of the latter. When we consider the relative Deligne cohomology, a relative class of this type is represented by a cocycle on $X$ together with a strong topological trivialization of its restriction to $A$ \cite{FR3}. Here we consider the trivialization as a suitable homotopy of differential functions.
\begin{Def} Given a differential function $(f, h, \omega): X \rightarrow (E_{n}, \iota_{n})$, a \emph{strong topological trivialization} of $(f, h, \omega)$ is a homotopy $(F, H, \pi^{*}\omega): X \times I \rightarrow (E_{n}, \iota_{n})$ between $(f, h, \omega)$ and a function of the form $(c_{e_{n}}, \chi(\rho), d\rho)$.
\end{Def}
It follows from the definition that $\omega = d\rho$. Analogously to definition \ref{GroupILinha}, the group $\hat{h}^{n}_{II'}(X, A)$ contains the homotopy classes of differential functions $(f, h, \omega): \Cyl(X, A) \rightarrow (E_{n}, \iota_{n})$ such that $(f, h, \omega)\vert_{\Cyl A}$ is a strong topological trivialization of $(f, h, \omega)\vert_{A}$. A homotopy $(F, H, \pi^{*}\omega): \Cyl(X, A) \times I \rightarrow (E_{n}, \iota_{n})$ between two such functions, whose trivialization on $A$ is $(c_{e_{n}}, \chi(\rho), d\rho)$, is required to satisfy $(F, H, \pi^{*}\omega)\vert_{A \times \{1\} \times I} = (c_{e_{n}}, \chi(\pi_{A}^{*}\rho), \pi_{A}^{*}d\rho)$. It follows that $\omega\vert_{A} = d\rho$. We can show that $\hat{h}^{n}_{II'}(X, A) \simeq \hat{h}_{II}^{n}(X, A)$ canonically, as we did for type I. In fact, there is a natural morphism $\varphi: \hat{h}_{II}^{n}(X, A) \rightarrow \hat{h}^{n}_{II'}(X, A)$ defined as follows: given $[(f, h, \omega)] \in \hat{h}_{II}^{n}(X, A)$, since by definition $(f, h, \omega)\vert_{A} = (c_{e_{n}}, \chi(\rho), d\rho)$, we extend $(f, h, \omega)$ to $\Cyl(X, A)$ in such a way that $(f, h, \omega)\vert_{\Cyl A} = (c_{e_{n}}, \chi(\pi_{A}^{*}\rho), \pi_{A}^{*}d\rho)$ (v.\ \cite{Upmeier} lemma 2.2). Such a morphism is actually an isomorphism: about the injectivity the proof of theorem \ref{IsoI} applies without variations, since, if $\varphi([(f, h, \omega)]) = 0$, it must hold that $\rho = 0$; about the surjectivity, the proof of theorem \ref{IsoI} can be adapted considering the homotopy $(F, H, \pi_{\Cyl\,A}^{*}\pi_{A}^{*}d\rho): A \times I \times I \rightarrow (E_{n}, \iota_{n})$ instead of $(F, H, 0)$.\footnote{Moreover, the definition of $K''$ becomes $K'' := \pi_{A \times \{1\} \times I}^{*}(H'' - \chi(\rho))$.} A remark analogous to the one after theorem \ref{IsoI} still holds. We can now construct the Bockstein map of \eqref{LongExact2}. We start from $[(f, h, \rho)] \in \hat{h}^{n-1}(A)$, and, by analogy with the Deligne cohomology, we must get a class in $\hat{h}_{II}^{n}(X, A)$ whose restriction to $A$ is $(c_{e_{n}}, \chi(-\rho), 0)$. Nevertheless, we cannot simply extend such a restriction to $X$, because the extension is not unique. Actually, if the Bockstein map were defined in this way, only the curvature $\rho$ of $[(f, h, \rho)]$ would be meaningful, therefore the kernel would be $\hat{h}_{\fl}^{n-1}(A)$, not the image of $\hat{h}_{\fl}^{n-1}(X)$. We thus need a different construction, passing through the group $\hat{h}_{II'}^{n}(X, A)$.

Briefly the idea is the following. We call $\pi_{1}: S^{1} \times A \rightarrow A$ the projection and we fix a marked point on $S^{1}$, e.g.\ $1 \in S^{1} \subset \mathbb{C}$. Given $[(f, h, \rho)] \in \hat{h}^{n-1}(A)$, we consider the unique class $[(F, H, dt \wedge \pi_{1}^{*}\rho)] \in \hat{h}_{I}^{n}(S^{1} \times A, \{1\} \times A)$ whose integral over $S^{1}$ is $[(f, h, \rho)]$, and, identifying the pair $(S^{1} \times A, \{1\} \times A)$ with the pair $(I \times A, \{0, 1\} \times A)$, we define $\beta^{n-1}[(f, h, \rho)] := [(F, H - \chi(t \cdot \pi_{A}^{*}\rho), 0)] \in \hat{h}_{II'}^{n}(X, A)$, the extension to $X$ being the trivial one. We describe this construction in more detail. 
\begin{itemize}
	\item Given $[(f, h, \rho)] \in \hat{h}^{n-1}(A)$, thanks to \cite[Lemma 4.4]{Upmeier} there exists a class $[(F, H, dt \wedge \pi_{1}^{*}\rho)] \in \hat{h}_{I}^{n}(S^{1} \times A, \{1\} \times A)$ such that $\int_{S^{1}} [(F, H, dt \wedge \pi_{1}^{*}\rho)] = [(f, h, \rho)]$. Such a class is unique: if we choose another one, the difference is a flat class $[(F', H', 0)] \in \hat{h}_{\fl}^{n}(S^{1} \times A, \{1\} \times A) \simeq \hat{h}_{\fl}^{n-1}(A)$, the isomorphism being given by the $S^{1}$-integration.\footnote{In particular, calling $A_{+} := A \sqcup \{\infty\}$ and $S$ the suspension, we have that $\hat{h}_{\fl}^{n}(S^{1} \times A, \{1\} \times A) \simeq \tilde{\hat{h}}_{\fl}^{n}((S^{1} \times A)/(\{1\} \times A)) \simeq \tilde{\hat{h}}_{\fl}^{n}(S(A_{+})) \simeq \tilde{\hat{h}}_{\fl}^{n-1}(A_{+}) \simeq \hat{h}_{\fl}^{n-1}(A)$.} Since $\int_{S^{1}} [(F', H', 0)] = 0$, we get $[(F', H', 0)] = 0$.
	\item Composing with the pull-back via the projection $p: (I \times A, \{0,1\} \times A) \rightarrow (S^{1} \times A, \{1\} \times A)$, we get a class represented by a differential function on $I \times A$ that we still call $(F, H, dt \wedge \pi_{A}^{*}\rho)$, whose restriction to $\{0,1\} \times A$ is $0$.
	\item We define the following differential function on $I \times A$:
\begin{equation}\label{ClassBock}
	(F, H - \chi(t \cdot \pi_{A}^{*}\rho), 0).
\end{equation}
Such a function is well-defined, since, by construction, $\delta H = \chi(dt \wedge \pi_{A}^{*}\rho) - F^{*}\iota_{n} = \chi(d(t \cdot \pi_{A}^{*}\rho)) - F^{*}\iota_{n} = \delta\chi(t \cdot \pi_{A}^{*}\rho) - F^{*}\iota_{n}$, hence $\delta(H - \chi(t \cdot \pi_{A}^{*}\rho)) = - F^{*}\iota_{n}$.
	\item By construction $(F, H - \chi(t \cdot \pi_{A}^{*}\rho), 0)\vert_{\{1\} \times A} = (c_{e_{n}}, \chi(-\rho), 0)$ (we recall that $d\rho = 0$ since $\rho$ is a curvature), and $(F, H - \chi(t \cdot \pi_{A}^{*}\rho), 0)\vert_{\{0\} \times A} = (c_{e_{n}}, 0, 0)$. Hence we extend such a class to $\Cyl(X,A)$ requiring that it is trivial on $X$, and we get a representative of a class in $\hat{h}_{II'}^{n}(X, A)$.
\end{itemize}
We can now verify that the map is well-defined and it is a group homomorphism. We have already pointed out that the class $[(F, H, dt \wedge \pi_{1}^{*}\rho)] \in \hat{h}_{I}^{n}(S^{1} \times A, \{1\} \times A)$ is unique. If we choose another representative $(F', H', dt \wedge \pi_{1}^{*}\rho)$, by definition there exists a homotopy $(\Phi, \mathcal{H}, dt \wedge \pi^{*}\rho)$ defined on $I \times S^{1} \times A$ which is constant on $\{1\} \times A$, therefore, pulling back to $I \times A$, we get a homotopy defined on $I \times I \times A$ which is constant on $\{0,1\} \times A$. Considering $(\Phi, \mathcal{H} - \chi(t \cdot \pi^{*}_{A}\rho), 0)$, we get a homotopy defined on $I \times I \times A$ which is constant and equal to $(c_{e_{n}}, 0, 0)$ on $\{0\} \times A$, and constant and equal to $(c_{e_{n}}, \chi(-\rho), 0)$ on $\{1\} \times A$. Therefore the two classes $(F, H - \chi(t \cdot \pi_{A}^{*}\rho), 0)$ and $(F', H' - \chi(t \cdot \pi_{A}^{*}\rho), 0)$ define the same class in $\hat{h}_{II'}^{n}(X, A)$. In order to show that it is a group homomorphism, we notice that the class represented by \eqref{ClassBock} is equal to $[(F, H, dt \wedge \pi_{A}^{*}\rho)] - [(c_{e_{n}}, \chi(t \cdot \pi_{A}^{*}\rho), dt \wedge \pi_{1}^{*}\rho)]$, therefore the Bockstein map we constructed is equal to the difference of two homomorphisms, hence it is a homomorphism. \\

\paragraph{}We can now prove the exactness of \eqref{LongExact2}:
\begin{itemize}
	\item \emph{Exactness in $\hat{h}^{n-1}(A)$.} Given $[(f, h, \rho)] \in \hat{h}^{n-1}(A)$, if $\beta^{n-1}[(f, h, \rho)] = 0$ then $\rho = 0$, since the restriction of \eqref{ClassBock} to $\{1\} \times A$ is $(c_{e_{n}}, -\chi(\rho), 0)$. Therefore the kernel of $\beta^{n-1}$ is contained in the flat part $\hat{h}_{\fl}^{n-1}(A)$, hence the exactness follows from the one of \eqref{LongExact1}.
	\item \emph{Exactness in $\hat{h}^{n}_{II}(X, A)$.} If a class belongs to the image of the Bockstein map, it follows from the previous construction that its restriction to $X$ is trivial. Viceversa, let us consider $[(F, H', \omega)] \in \hat{h}^{n}_{II'}(X, A)$ such that $[(F, H', \omega)]\vert_{X}$ is trivial. By the homotopy extension property applied to the pair $(\Cyl(X,A), X \cup (\{1\} \times A))$, we can suppose that $(F, H', \omega)\vert_{X} = (c_{e_{n}}, 0, 0)$. By definition of $\hat{h}^{n}_{II'}(X, A)$, we have that $(F, H', \omega)\vert_{\{1\} \times A} = (c_{e_{n}}, -\chi(\rho), 0)$, therefore we define $H := H' + \chi(t \cdot \pi_{A}^{*}\rho)$, so that $(F, H', \omega)$ takes the form \eqref{ClassBock}. It follows that the differential function $(F, H, dt \wedge \pi_{A}^{*}\rho)$ is well-defined on $I \times A$, since $\delta H = \delta H' + \chi(dt \wedge \pi_{A}^{*}\rho) = -F^{*}\iota_{n} + \chi(dt \wedge \pi_{A}^{*}\rho)$. Being $(F, H, dt \wedge \pi_{A}^{*}\rho)\vert_{\partial I \times A} = (c_{e_{n}}, 0, 0)$, it represents a class in $\hat{h}^{n}_{I}(I \times A, \partial I \times A) \simeq \hat{h}^{n}_{I}(S^{1} \times A, \{1\} \times A)$, hence we obtain $[(F, H, dt \wedge \pi_{1}^{*}\rho)] \in \hat{h}^{n}_{I}(S^{1} \times A, \{1\} \times A)$. Integrating it over $S^{1}$ we get a class $[(f, h, \rho)] \in \hat{h}^{n-1}(A)$ such that $\beta^{n-1}[(f, h, \rho)] = [(F, H', \omega)]$.\footnote{Of course the class $[(f, h, \rho)]$ depends in general on the representative $(F, H', \omega)$ chosen, since, when passing to $S^{1} \times A$ from $I \times A$, there is not a well-defined push-forward in cohomology. This is not a problem, since we must find at least one class in $\hat{h}^{n-1}(A)$ whose image under the Bokstein map is $[(F, H', \omega)]$. Actually, it cannot be unique in general because of exactness of \eqref{LongExact2}.}
	\item \emph{Exactness in $\hat{h}^{n}(X)$.} By construction the restriction to $A$ of a class in $\hat{h}^{n}_{II}(X, A)$ is topologically trivial. Viceversa, given a class $[(f, h, \omega)] \in \hat{h}^{n}(X)$ which is topologically trivial on $A$, by the homotopy extension property we can suppose that $(f, h, \omega)\vert_{A} = (c_{e_{n}}, \chi(\rho), d\rho)$. It follows that such a class is the image of a class $[(f, h, \omega)] \in \hat{h}^{n}_{II}(X, A)$.
\end{itemize}

\subsection{Type III} We briefly consider the groups of type III in order to complete the picture. When we consider the relative Deligne cohomology, a relative class of this type is represented by a cocycle on $X$ together with a topological trivialization of its restriction to $A$. Here we consider the trivialization as a suitable homotopy of differential functions.
\begin{Def} Given a differential function $(f, h, \omega): X \rightarrow (E_{n}, \iota_{n})$ a \emph{topological trivialization} of $(f, h, \omega)$ is a homotopy $F: X \times I \rightarrow E_{n}$ between $f$ and $c_{e_{n}}$.
\end{Def}
The group $\hat{h}^{n}_{III'}(X, A)$ contains the homotopy classes of differential functions $(f, h, \omega): X \rightarrow (E_{n}, \iota_{n})$ with a topological trivialization of $f\vert_{A}$ extending $f$ to $\Cyl(X, A)$. A homotopy $(F, H, \pi^{*}\omega): X \times I \rightarrow (E_{n}, \iota_{n})$ between two such functions must extend to a homotopy of homotopies $F$ on $\Cyl \Cyl A$ such that $F\vert_{A \times \{1\} \times I} = c_{e_{n}}$. As in the previous cases, there is a natural morphism $\varphi: \hat{h}_{III}^{n}(X, A) \rightarrow \hat{h}^{n}_{III'}(X, A)$ which is actually an isomorphism. In order to construct the Bockstein map of \eqref{LongExact3}, we consider a class $[f] \in h^{n}(A)$ and the corresponding class $[F] \in h^{n+1}(SA)$. We pull-back such a class to $[F] \in h^{n+1}(I \times A, \partial I \times A)$, and we extend $F$ to a trivial differential function on $X$. In this way we get a class in $\hat{h}^{n}_{III'}(X, A)$ whose restriction to $X$ is trivial. One can verify that \eqref{LongExact3} is exact.

\subsection{Type IV} Let us consider the group $\hat{h}^{n}_{II}(X, A)$: we define the subgroup $\hat{h}^{n}_{II,\ch}(X, A)$ made by those classes $[(f, h, \omega)]$ such that $[(f, h, \omega)]\vert_{A} = 0$: by definition this means that $(f, h, \omega)\vert_{A} = (c_{e_{n}}, \chi(\rho), 0)$, with $\rho$ a closed form representing a class belonging to the image of the Chern character. This subgroup is bigger than $\hat{h}^{n}_{I}(X, A)$, because in the latter the form $\rho$ must be $0$ (this means that, in $\hat{h}^{n}_{I}(X, A)$, the function $(f, h, \rho)$ must vanish on $A$ as a single representative, not only as a cohomology class). There is a natural morphism:
\begin{equation}\label{MorphIIChIV}
	\varphi: \hat{h}^{n}_{II, \ch}(X, A) \rightarrow \hat{h}^{n}_{IV}(X, A),
\end{equation}
sending $[(f, h, \omega)]$ to $[(f, [h], \omega)]$. The latter is well-defined, since a homotopy of type II representatives restricts to $(c_{e_{n}}, \chi(\pi_{A}^{*}\rho), 0)$ on $A \times I$, hence, being $(c_{e_{n}}, [\chi(\pi_{A}^{*}\rho)], 0) = (c_{e_{n}}, [0], 0)$, it also defines a homotopy of type IV representatives.

We denote by $\Omega^{n}_{\ch}(X, \mathfrak{h}^{\bullet}_{\mathbb{R}})$ the group of $\mathfrak{h}^{\bullet}_{\mathbb{R}}$-valued closed forms of degree $n$ that represent a class belonging to the image of the Chern character. There is a natural map:
\begin{equation}\label{MorphPsi}
	\psi: \Omega^{n-1}_{\ch}(X, \mathfrak{h}^{\bullet}_{\mathbb{R}}) \rightarrow \hat{h}^{n}_{II, \ch}(X, A)
\end{equation}
defined by $\psi(\rho) = [(c_{e_{n}}, \chi(\rho), 0)]$.
\begin{Lemma}\label{MorphIIChIVThm} The morphism \eqref{MorphIIChIV} is surjective and its kernel is the image of \eqref{MorphPsi}.
\end{Lemma}
\paragraph{Proof:} For the surjectivity, given a class $[(f, [h], \omega)] \in \hat{h}^{n}_{IV}(X, A)$, we choose a representative $(f, h, \omega)$. Then $h\vert_{A} \in Z^{n-1}_{\ch}(A, \mathbb{R})$ and it is cohomologous to a closed differential form, i.e.\ $h\vert_{A} = \chi(\rho) + \delta k$. The function $(c_{e_{n}}, h\vert_{A}, 0)$ is homotopic to $(c_{e_{n}}, \chi(\rho), 0)$ via a homotopy of the form $(c_{e_{n}}, \chi(\rho) + \delta K, 0)$ (v.\ \cite[Lemma 2.3]{Upmeier}). By the homotopy extension property, we get a type II representative whose class is sent to $[(f, [h], \omega)]$ by $\varphi$.

The inclusion $\IIm \psi \subset \Ker \varphi$ is an immediate consequence of the definition of $\hat{h}^{n}_{IV}(X, A)$. Viceversa, let us suppose that $\varphi[(f, h, 0)] = 0$. Then there is a type IV homotopy $(F, H, 0)$ between $(f, h, 0)$ and $(c_{e_{n}}, \chi(\rho'), 0)$, where $\chi(\rho') \in Z^{n}_{\ch}(X, \mathbb{R})$. We have that $h\vert_{A} = (c_{e_{n}}, \chi(\rho), 0)$ and, since $A$ is closed in $X$, we can suppose that $\rho'\vert_{A} = \rho$. By construction $F\vert_{A \times I} = c_{e_{n}}$. Moreover, $H\vert_{A \times I} \in Z_{\ch}^{n-1}(A \times I, \mathbb{R})$ and, being $A$ a deformation retract of $A \times I$, we have $H\vert_{A \times I} = \pi^{*}\chi(\rho) + \delta K$. We extend $K$ to the whole $X \times I$ (imposing that it gives $0$ when evaluated on a simplex not contained in $A \times I$), and we consider the homotopy $(F, H - \delta K, 0)$ between $(f, h, 0)$ and $(c_{e_{n}}, \chi(\rho'), 0)$, which is constant on $A \times I$. This shows that $[(f, h, 0)] = \psi(\rho')$. $\square$ \\

In order to prove the exactness of \eqref{LongExact4}, we generalize to any cohomology theory the groups defined by formula (10) of \cite{FR3}, i.e.:
\begin{equation}\label{hBar}
	\overline{h}^{n}(X, A) := \frac{\Ker(\hat{h}^{n}_{II}(X, A) \rightarrow \hat{h}^{n}(A))}{\IIm(\hat{h}^{n-1}(X) \rightarrow \hat{h}^{n}_{II}(X, A))}.
\end{equation}
\begin{Theorem}\label{ThmXip3} There is a natural isomorphism:
\begin{equation}\label{Xip3}
	\Xi^{n}: \overline{h}^{n}(X, A) \overset{\simeq}\longrightarrow \hat{h}^{n}_{IV}(X, A).
\end{equation}
\end{Theorem}
\paragraph{Proof:} The numerator of \eqref{hBar} corresponds to $\hat{h}^{n}_{II, \ch}(X, A)$, because the latter is exactly the group of classes of type II vanishing on $A$. Let us show that the denominator coincides with the image of \eqref{MorphPsi}, so that the result follows from lemma \ref{MorphIIChIVThm}. Because of \eqref{LongExact2}, the kernel of the map $\hat{h}^{n-1}(A) \rightarrow \hat{h}^{n}_{II}(X, A)$ is made by the image of flat classes on $X$, therefore, when composing with $\hat{h}^{n-1}(X) \rightarrow \hat{h}^{n-1}(A)$, only the curvature of the original class in $X$ is meaningful. In particular, given $[(f, h, \rho)] \in \hat{h}^{n-1}(A)$, in order to compute the Bockstein map of \eqref{LongExact2}, we consider the unique class $[(F, H, dt \wedge \pi_{1}^{*}\rho)] \in \hat{h}_{I}^{n}(S^{1} \times A, \{1\} \times A)$ whose integral over $S^{1}$ is $[(f, h, \rho)]$, and we define $\beta^{n-1}[(f, h, \rho)] := [(F, H - \chi(t \cdot \pi_{A}^{*}\rho), 0)] \in \hat{h}_{II'}^{n}(X, A)$. If $[(f, h, \rho)]$ is the restriction of a class on the whole $X$, we apply the same procedure to the whole class, obtaining a homotopy $(F, H - \chi(t \cdot \pi_{X}^{*}\rho), 0)$ on $I \times X$. The latter restricts on $\Cyl(X,A)$ to a representative of $\beta^{n-1}[(f, h, \rho)]$, and restricts on $\{1\} \times X$ to $(c_{e_{n}}, -\chi(\rho), 0)$. It follows from the remark after theorem \ref{IsoI} (adapted to type II) that $\beta^{n-1}[(f, h, \rho)] = \psi(-\rho)$. $\square$ \\

The fact that \eqref{LongExact4} is exact therefore follows from \cite[Theorem 2.1]{FR3}, which we repeat here for completeness with the notation of the present paper.
\begin{Theorem}\label{ExactnessIV} The sequence \eqref{LongExact4} is exact.
\end{Theorem}
\paragraph{Proof:} The Bockstein map $\beta$ of \eqref{LongExact4} is induced from the one of \eqref{LongExact2}, that we call $\beta'$: the image of $\beta'$ is contained in the numerator of \eqref{hBar} because of the exactness of \eqref{LongExact2}. By definition of the denominator of \eqref{hBar}, the kernel of $\beta$ is the image of the restriction map $\hat{h}^{p-1}(X) \rightarrow \hat{h}^{p-1}(A)$, thus \eqref{LongExact4} is exact in $\hat{h}^{p-1}(A)$. The image of $\beta'$ contains the denominator of \eqref{hBar}, therefore the exactness in $\hat{h}^{p}_{IV}(X, A)$ follows from the one of \eqref{LongExact2}. Finally, in order to prove the exactness in $\hat{h}^{p}(X)$, we consider the following commutative diagram:
	\[\xymatrix{
	\hat{h}^{p}_{I}(X, A) \ar[r]^(.55){\eta} \ar[d]_{\psi} & \hat{h}^{p}(X) \\
	\hat{h}^{p}_{IV}(X, A), \ar[ur]_(.55){\nu}
}\]
where $\eta$ is the map appearing in \eqref{LongExact1}, $\nu$ the one appearing in \eqref{LongExact4} and $\psi$ is the composition of the embedding $\hat{h}^{p}_{I}(X, A) \rightarrow \hat{h}^{p}_{II}(X, A)$, whose image is contained in the numerator of \eqref{hBar} because of the exactness of \eqref{LongExact1}, with the projection to the quotient in \eqref{hBar}. We show that $\IIm \, \eta = \IIm \, \nu$, so that the exactness of \eqref{LongExact4} in $\hat{h}^{p}(X)$ follows from the one of \eqref{LongExact1}. Obviously $\IIm \, \eta \subset \IIm \, \nu$. For the converse, the image of the embedding $\hat{h}^{p}_{I}(X, A) \rightarrow \hat{h}^{p}_{II}(X, A)$ is the subset of classes which are trivial when pulled-back to $\hat{h}^{p}(A)$, therefore, applying $\nu$ to the numerator of \eqref{hBar}, we get classes belonging to the kernel of $\hat{h}^{p}(X) \rightarrow \hat{h}^{p}(A)$, i.e.\ to the image of $\eta$. $\square$

\section{Generic map}\label{GenericMap}

Up to now we have considered the relative groups of a pair $(X, A)$, the latter being equivalent to an embedding $i: A \hookrightarrow X$ of a closed submanifold. We are going to generalize the definition to any smooth map of manifolds $\varphi: A \rightarrow X$. The most obvious generalizations are the following:
\begin{itemize}
	\item the cylinder $\Cyl(X, A)$ is replaced by $\Cyl(\varphi)$, the latter being the union of the two spaces $X$ and $A \times I$, identifying $(a, 0) \in A \times I$ with $\varphi(a) \in X$;
	\item a restriction of a class on $X$ to $A$ is replaced by the pull-back via $\varphi$.
\end{itemize}
Nevertheless, we cannot consider the direct generalization of the groups $\hat{h}_{I}^{n}(X, A)$, $\ldots$, $\hat{h}_{IV}^{n}(X, A)$, because there is not a suitable version of lemma \ref{HEP} if $\varphi$ is not a cofibration. The only possibility consists of considering the groups $\hat{h}_{I'}^{n}(X, A)$, $\ldots$, $\hat{h}_{III'}^{n}(X, A)$, defined in section \ref{ExSequences}, as the original definition of the relative groups. Such definitions are similar to the one of topological relative cohomology via the cone of $\varphi$.

In the following, in order to define differential functions $(f, h, \omega): \Cyl(\varphi) \rightarrow (E_{n}, \iota_{n})$, we argue similarly to the case of a closed embedding (v.\ comments at the beginning of section \ref{HomExtProp}). In particular, we replace the embedding $\iota: \Cyl(X, A) \hookrightarrow X \times I$ with the map $\iota: \Cyl(\varphi) \rightarrow X \times I$, defined by $\iota([x]) := (x, 0)$ and $\iota[(a, t)] := (\varphi(a), t)$. Moreover, we consider the maps:
	\[\iota_{\Cyl\,A}: \Cyl\,A \rightarrow \Cyl(\varphi) \qquad \iota_{A}: A \rightarrow \Cyl(\varphi),
\]
defined by $\iota_{\Cyl\,A}(a, t) = [(a, t)]$ and $\iota_{A}(a) = [(a, 0)] = [\varphi(a)]$.

\subsection{Type I} The group $\hat{h}^{n}_{I'}(\varphi)$ contains the homotopy classes of differential functions $(f, h, \omega): \Cyl(\varphi) \rightarrow (E_{n}, \iota_{n})$ such that $\iota_{\Cyl\,A}^{*}(f, h, \omega)$ is a geometric trivialization of $\iota_{A}^{*}(f, h, \omega)$. It follows that $\varphi^{*}\omega = 0$. A homotopy $(F, H, \pi^{*}\omega): \Cyl(\varphi) \times I \rightarrow (E_{n}, \iota_{n})$ between two such functions is required to satisfy $(F, H, \pi^{*}\omega)\vert_{A \times \{1\} \times I} = (c_{e_{n}}, 0, 0)$. The corresponding long exact sequence is analogous to \eqref{LongExact1}, the Bockstein map being constructed as stated at the end of subsection \ref{TypeISubS}. The only difference is that we have to extend the class we get on $\Cyl\,A$, which is trivial on $A \times \{0\}$, to $\Cyl(\varphi)$ instead of $\Cyl(X, A)$. The exactness in $\hat{h}^{n-1}(A)$ follows from the long exact sequence of the flat theory. In $\hat{h}^{n}_{I}(\varphi)$ it is straightforward from the construction of the Bockstein map (if the restriction of a class to $X$ is zero, we get a flat class on $\Cyl(A)$ which is trivial on the two bases, hence it is the image, via the suspension isomorphism, of a flat class on $A$). Finally, exactness in $\hat{h}^{n}(X)$ can be proven in the following way. Given $\alpha \in \hat{h}^{n}(X)$, if $\varphi^{*}\alpha = 0$, we can glue in $\Cyl(\varphi)$ a differential function on $X$ representing $\alpha$ and a homotopy between $\varphi^{*}\alpha$ and $(c_{e_{n}}, 0, 0)$. We thus get a class in $\hat{h}^{n}_{I}(\varphi)$, whose restriction to $X$ is $\alpha$.

\subsection{Type II} The group $\hat{h}^{n}_{II'}(\varphi)$ contains the homotopy classes of differential functions $(f, h, \omega): \Cyl(\varphi) \rightarrow (E_{n}, \iota_{n})$ such that $\iota_{\Cyl\,A}^{*}(f, h, \omega)$ is a strong topological trivialization of $\iota_{A}^{*}(f, h, \omega)$. It follows that $\varphi^{*}\omega = d\rho$. A homotopy $(F, H, \pi^{*}\omega): \Cyl(\varphi) \times I \rightarrow (E_{n}, \iota_{n})$ between two such functions, whose trivialization on $A$ is $(c_{e_{n}}, \chi(\rho), d\rho)$, is required to satisfy $(F, H, \pi^{*}\omega)\vert_{A \times \{1\} \times I} = (c_{e_{n}}, \chi(\pi_{A}^{*}\rho), \pi_{A}^{*}d\rho)$. The corresponding long exact sequence is analogous to \eqref{LongExact2}, the Bockstein map being constructed as in the case of a closed embedding. The homotopy extension property is used only in the proof of the exactness in $\hat{h}^{n}_{II'}(\varphi)$, with respect to the pair $(\Cyl(\varphi), X \cup (\{1\} \times A))$. In this case we can apply a suitable generalization of lemma \ref{HEP}, since $X \cup (\{1\} \times A)$ is a closed subspace of $\Cyl(\varphi)$. In particular, we call ``manifold with generalized cylinders'' a space defined inductively in the following way:
\begin{itemize}
	\item we start from a manifold $X_{0}$; we define $\iota_{0}: X_{0} \rightarrow X_{0}$ as the identity map.
	\item We are given inductively a space $X_{n}$ and a map $\iota_{n}: X_{n} \rightarrow X_{0} \times I^{n}$. We consider a manifold $A_{n}$ and a smooth map $\varphi_{n}: A_{n} \rightarrow X_{n}$, where ``smooth'' means that $\iota_{n} \circ \varphi_{n}: A_{n} \rightarrow X_{0} \times I^{n}$ is smooth. We define $X_{n+1} := \Cyl(\varphi_{n})$. Moreover, $\iota_{n+1}: X_{n+1} \rightarrow X_{0} \times I^{n+1}$ is defined as follows: $\iota_{n+1}([x]) = (\iota_{n}(x), 0)$ for $x \in X_{n}$, and $\iota_{n+1}[(a, t)] = (\iota_{n} \circ \varphi_{n}(a), t)$.
\end{itemize}
After a finite number of steps we get $X = X_{n}$ and $\iota = \iota_{n}: X \rightarrow X_{0} \times I^{n}$. Moreover, if $Y$ is a manifold, $f: Y \rightarrow X$ is smooth if and only if $\iota \circ f$ is. When also $Y$ is a manifold with generalized cylinders, $f: Y \rightarrow X$ is smooth if and only if, for any manifold $Z$ and any smooth map $\xi: Z \rightarrow Y$, the composition $f \circ \xi$ is smooth.

If $(X, A)$ is a pair of manifolds with generalized cylinders (hence $A \subset X$ is a closed embedding, even if $X$ and $A$ may contain generalized cylinders), lemma \ref{HEP} holds with the same proof. Therefore, also the exactness of \eqref{LongExact2} holds (in order to prove the exactness in $\hat{h}^{n}(X)$, we apply the same argument we used about type I).

\subsection{Type III} The definition of the group $\hat{h}^{n}_{III'}(\varphi)$ is analogous to the previous ones.

\subsection{Type IV} For type IV, the immediate generalization of lemma \ref{MorphIIChIVThm} does not hold. In particular, the proof of the fact that $\Ker \varphi \subset \IIm \psi$ fails, since, given a form $\rho$ on $A$, in general we cannot find a form $\rho'$ on $X$ such that $\varphi^{*}\rho' = \rho$ (for example, when $A \subset X$ is not closed, in general we cannot extend a form from $A$ to $X$). Hence, in the proof of the lemma, when we find a homotopy between $(f, h, 0)$ and $(c_{e_{n}}, \chi(\rho'), 0)$, we cannot guarantee that $\varphi^{*}\rho' = \rho$, thus we cannot find a type II homotopy between $(f, h, 0)$ and $\psi(\rho')$. The only possibility seems to consist of generalizing formula \eqref{hBar}, defining $h^{n}_{IV}(\varphi) := \Ker(\hat{h}^{n}_{II}(\varphi) \rightarrow \hat{h}^{n}(A))/\IIm(\hat{h}^{n-1}(X) \rightarrow \hat{h}^{n}_{II}(\varphi))$. Theorem \ref{ExactnessIV} holds, since the proof only relies on the exactness of the type I and type II sequences.


\end{document}